\documentclass[12pt]{amsart}
\usepackage{amssymb}
\usepackage{amssymb,amsmath,amscd}

\def\+{{+\!\!\!+}}

\def\d{\partial}

\def\g{\gamma}
\def\G{\Gamma}

\def\D{{\cal D}}

\def\pmb#1{\setbox0=\hbox{#1}%
\kern.0em\copy0\kern-\wd0
\kern-.04em\copy0\kern-\wd0
\kern.08em\copy0\kern-\wd0
\kern-.04em\raise.0433em\box0 }         %poor man's bold macro (TexBook)
\def\rank{\textstyle{\rm{rank}}}

\def\id{{\rm id}}
\def\Ad{{\rm Ad}}

\def\bis{{\rm Bis}}

\def\X{\hat{X}}
\def\A{{\mathcal A}}

\newcommand{\nc}{\newcommand}
\nc{\beq}{\begin{equation}}
\nc{\eeq}[1]{\label{#1}\end{equation}}
\nc{\ber}{\begin{eqnarray}}
\nc{\eer}[1]{\label{#1}\end{eqnarray}}
\nc{\pek}[1]{\cite{#1}}
\nc{\enr}[1]{(\ref{#1})}
\nc{\kal}[1]{{\cal{#1}}}
\nc{\dott}{\;\cdot\;}

\def\G{{\mathcal G}}
\def\M{{\mathcal M}}
\def\H{{\mathcal H}}
\def\g{{\mathbf g}}
\def\h{{\mathbf h}}

\def\sphere{{\mathbf S}}
\def\D{{\mathcal D}}
\def\E{{\mathcal E}}
\def\Ep{{E'{}}}
\def\endproof{\hskip5pt \vrule height4pt width4pt depth0pt \par}
\newcommand{\be}{\begin{equation}}
\newcommand{\ee}{\end{equation}}
\newcommand{\bea}{\begin{eqnarray}}
\newcommand{\eea}{\end{eqnarray}}

\theoremstyle{remark}
\newtheorem*{Ack}{Acknowledgment}

\theoremstyle{definition}
\newtheorem{Thm}{Theorem}[section]

\newtheorem{Lem}[Thm]{Lemma}

\newtheorem*{Thm*}{Theorem}

\newtheorem{Def}[Thm]{Definition}
\newtheorem{Exa}[Thm]{Example}

\newtheorem{pro}[Thm]{Proposition}

\newtheorem{rmk}[Thm]{Remark}
\newtheorem{defi}[Thm]{Definition}

\begin{document}

\title[Lie algebroids, Lie groupoids and TFT]
{Lie algebroids, Lie groupoids and TFT}
\date{2005}

\author[F.~Bonechi]{Francesco Bonechi}
\author[M.~Zabzine]{Maxim Zabzine}
\address{I.N.F.N. and Dipartimento di Fisica,
  Via G. Sansone 1, 50019 Sesto Fiorentino - Firenze, Italy}
\email{Francesco.Bonechi@fi.infn.it}
\address{Department of Theoretical Physics,
 Uppsala University, Box 803, SE-75108 Uppsala, Sweden}
 \address{Kavli Institute of Theoretical Physics, University of California, Santa Barbara, CA 93106 USA}
\email{Maxim.Zabzine@teorfys.uu.se}

\date{December 2005, Preprint numbers: UUITP-22/05, NSF-KITP-05-112}
%\subjclass[2000]{}
%\keywords{}
\thanks{The research of M.~Z. was supported by VR Grant No. 621-2004-3177
 and in part by NSF Grant No. PHY99-07949.}

\begin{abstract}
 We construct the moduli spaces associated to the solutions of equations of
  motion (modulo gauge transformations) of the Poisson sigma model with target
   being an  integrable Poisson  manifold. The construction can be easily extended
   to a case of a generic integrable Lie algebroid.  Indeed for any Lie algebroid one can associate
   a BF-like topological field theory which localizes on the space of algebroid morphisms, that can be seen as a   generalization of flat connections to the groupoid case.
     We discuss the finite gauge transformations and discuss the corresponding moduli spaces.
    We consider the theories both without and with boundaries.
\end{abstract}

\maketitle

\section{Introduction}\label{sec-intro}

Topological field theory (TFT) plays  a prominent role in the
investigation of geometry and topology  of the moduli spaces of flat
connections over a two dimensional surface $\Sigma$.  In particular these moduli spaces
  appear as the phase space of Chern-Simons theory, the localization locus of two dimensional
Yang-Mills and as the stationary points of $BF$ theory. The
application of quantum field theoretical methods has produced new
results such as the formulas for the symplectic volume and for the
intersection numbers (e.g., see \cite{Blau:1993hj}).

The Poisson sigma model (PSM) is another example of two dimensional TFT
  introduced in \cite{Ikeda:1993fh,Schaller:1994es} which is a
sigma model defined on a two dimensional surface $\Sigma$ with
target being a Poisson manifold. The $BF$-theory
  and $A$-model are particular examples of PSM.  Recently PSM
 has attracted additional attention due to its relation with the
deformation quantization
\cite{Cattaneo:1999fm}. However the potential use of the PSM as a
TFT still remains to be investigated.

In this paper we show how a generalization of the moduli space of
flat connections over $\Sigma$ naturally appears when we study the
stationary configurations of PSM. We complete the study of
the moduli space of stationary points of PSM modulo the gauge
transformations initiated in \cite{bonechi} for the special
case of Poisson-Lie groups. In this situation the equations of motion
of the model have a straightforward geometrical interpretation as
the equations of a flat connection on the trivial bundle for the
dual group and a parallel section for the fibre bundle associated
to the dressing action. This observation dictates an obvious
choice for the on-shell finite gauge transformations and the
corresponding definition of the moduli space of solutions.
Moreover  the resulting space has a natural description where for every
symplectic leaf of the target we associate the moduli space of flat
connections for the isotropy group of the leaf.

In the present work we address the general case. Surprisingly it
turns out that once we properly identify the gauge transformations
the moduli space admits the same description. Furthermore
 the geometrical interpretation of the equations of PSM
as  algebroid morphisms between $T\Sigma$ and $T^*\M$
 suggests that  the discussion remains valid if we
substitute $T^*\M$ by a generic (integrable) Lie algebroid $E$.
This extension leads us beyond the PSM and puts the results
 in a different perspective.

An algebroid over a point is  a Lie algebra and the moduli
space of algebroid morphisms coincides with the moduli space of
flat connections for the trivial bundle. Thus in the general case we
can look at the equations of motion as a generalized flat
connection equation where the structure constants depend on the
base manifold. By applying the Lie theorems for integrable
algebroids, we can equivalently deal with groupoid morphisms from
the fundamental groupoid $\Pi(\Sigma)$ to the groupoid $\G(E)$
integrating $E$. This must be seen as the generalization of the
holonomy description of the moduli space of flat $G$-connections
as the space of $G$-representations of the fundamental group. Our
result describes this generalized moduli space as the union over
the leaves of the representations of the fundamental group in the
isotropy group.

The structure of the paper is as follows. In Section \ref{DEF} we review some basic notions
 from algebroid and groupoid theory.  In
Section \ref{algebroids_TFT} we discuss the equations of motion and the corresponding
 TFTs. Also we consider the natural boundary conditions in this setup.
Section \ref{GAUGE} is devoted to the finite
gauge transformations which form a groupoid.  We explain
 their relation to the algebroid (groupoid) homotopy.
  In Section \ref{MODULI} we
discuss the moduli spaces and their
 equivalence to the various generalizations of flat
connections.
 Section \ref{END} contains the summary and the list of open problems.

\begin{Ack}
We are grateful to Alberto Cattaneo, Carlo Rossi, Marco Zambon for
useful discussions and to Gabriele Vezzosi for pointing out to us
the considerations contained in Section \ref{END}.
\end{Ack}

\section{Lie algebroids and Lie groupoids}
\label{DEF}

In this Section we recall  some basic notions from the theory
of Lie algebroids and Lie groupoids. For a more extensive
discussion we refer to \cite{CW} and \cite{Mak}.

\subsection{Lie algebroids}\label{La}

\Def{ A {\it Lie algebroid} $(E, \M, \rho, \{\,\,,\,\,\})$ is a
vector bundle $E$ over a manifold
 $\M$ together with a bundle map (the anchor) $\rho: E \rightarrow T \M$ and a Lie bracket
 $\{\,\,,\,\,\}$ on the space $\Gamma(E)$ of sections of $E$ satisfying the compatibility condition
\begin{equation}\label{compalgerb}
 \{v,fu\} = f\{v,u\} + {\mathcal L}_{\rho_* v} fu,\,\,\,\,\,\,\,\,\,\,\,u,v \in \Gamma(E),\,\,\,\,f\in
 C^\infty({\M})
\end{equation}
 where $\rho_* : \Gamma(E) \rightarrow \Gamma(T\M)$ is the induced map of sections
  and ${\mathcal L}$  is the Lie derivative}.

 It follows from the definition that $\rho_*$ is a morphism of Lie
 algebras. On a trivializing chart $U$ we can choose the local coordinates
$X^\mu$ ($\mu =1,..., \dim {\M}$)
 and a basis $e^A$, ($A=1,...,\rank E$) on the fiber (e.g., the basis of constant sections on $E|_U$).
 In these local coordinates we introduce the anchor $\rho^{\mu A}$ and the structure functions
\begin{equation}\label{definacst}    
 \rho(e^A)(X) = \rho^{\mu A}(X) \d_\mu,\,\,\,\,\,\,\,\,\,\,\,\,\,\,\,
 \{ e^A, e^B\} = f^{AB}_{\,\,\,\,\,\,\,\,\,C} e^C.
\end{equation}
 The compatibility condition (\ref{compalgerb}) implies the following equations
\ber
 \label{consalg1}\rho^{\nu A} \d_\nu \rho^{\mu B} - \rho^{\nu B} \d_\nu \rho^{\mu A}=
  f^{AB}_{\,\,\,\,\,\,\,\,\,C} \rho^{\mu C}\\
\label{consalg2} \rho^{\mu [D} \d_\mu f^{AB]}_{\,\,\,\,\,\,\,\,\,C} +
 f^{[AB}_{\,\,\,\,\,\,\,\,\,L} f^{D]L}_{\,\,\,\,\,\,\,\,\,C}  =0 ,
\eer{localcoordinate}
 where $[\,\,\, ]$ stands for the antisymmetrization.

To any Lie algebroid we associate a characteristic foliation. The
isotropy algebra for $x\in \M$ is defined as the kernel of the
anchor map $\rho$
\begin{equation}\label{kernlanch}
 {\mathbf g}_x = \ker (\rho|_{E_x}).
\end{equation}
The characteristic foliation is a singular foliation on $\M$
determined by the distribution $x\rightarrow\mbox{Im}
(\rho|_{E_x})$. For any $x, y$ in the same leaf $\in L$ we have
${\mathbf g}_x \simeq {\mathbf g}_y$. Hence we get a bundle of Lie
algebras over $L$
\begin{equation}\label{budnlealg}
 {\mathbf g}_L = \bigcup\limits_{x \in L} {\mathbf g}_x\,\,\,\rightarrow\,\,\,L.
\end{equation}

Here are some examples of Lie algebroids which will be relevant
for further discussion.

\begin{Exa}\label{example_algbd1} Every Lie algebra $\mathbf g$ is
an algebroid over a point ($\rho=0$).
 \end{Exa}
\begin{Exa}\label{example_algbd2}
  The tangent bundle $T \M$ of a smooth manifold $\M$ is an
algebroid with the bracket between vector fields and $\rho=\id$.
 \end{Exa}
\begin{Exa}\label{example_algbd3}
 Let $\gamma:\mathbf g\rightarrow {\rm Vect}(\M)$ be a right
action of a Lie algebra $\mathbf g$ on a manifold $\M$. The {\it
action Lie algebroid} is defined on $\M\times{\mathbf g}$ over
$\M$ with anchor $\rho(m,\xi)=\gamma(\xi)(m)\in T_m\M$ and bracket
between $v,w\in \Gamma(\M\times\g)=C^\infty(\M,\g)$
$$\{v,w\}(m)=[v(m),w(m)] + {\mathcal
L}_{\gamma(v(m))}(w)(m)-{\mathcal L}_{\gamma(w(m))}(v)(m)\;,$$
where $[\ ,]$ denotes the bracket in $\mathbf g$ and $\mathcal L$
the Lie derivative.
 \end{Exa}
\begin{Exa}\label{example_algbd4}
  Let $\M$ be a Poisson manifold with Poisson tensor
$\alpha\in\Gamma(\wedge^2T\M)$. The associated canonical algebroid
is defined on $T^*\M$ by choosing as anchor the contraction
$\sharp_\alpha$ of cotangent vectors with $\alpha$ and by defining
the bracket on exact forms as $\{df,dg\}=d\{f,g\}, f,g\in
C^\infty(\M)$ and extending it to all $\Gamma(T^* \M)$ with
(\ref{compalgerb}).
 \end{Exa}

Next following Higgins and Mackenzie \cite{higgins}  we give a definition of Lie algebroid morphism
 which plays a central role in our investigation:

\begin{defi} \label{def_algbd_morph} Let $(E_1, \M_1, \rho_1, \{\,\,,\,\,\}_1)$ and $(E_2,
\M_2, \rho_2, \{\,\,,\,\,\}_2)$
 be Lie algebroids. Then {\it a morphism of Lie algebroids} is a vector bundle morphism
\begin{equation}
\begin{CD}
 E_1 @>\Phi >>   E_2     \\
@V \pi_1  VV   @V \pi_2 VV \\
 \M_1 @>\phi>> \M_2
\end{CD}
\end{equation}
 such that
\begin{equation}\label{commancmap} \rho_2 \circ \Phi  = d\phi\circ \rho_1,
\end{equation}
 where $d\phi: T\M_1 \rightarrow T\M_2$ and such that for
 arbitrary $V, W \in \Gamma(E_1)$ with $\Phi$-decomposition
\begin{equation}\label{decompof}
 \Phi\circ V = \sum V^i (e_i\circ\phi),\,\,\,\,\,\,\,\,\,\,\,\,\,\,\,\,\,\,
 \Phi \circ W = \sum W^i (\tilde{e}_i \circ\phi)
\end{equation}
where $W^i, V^i \in C^\infty(\M_1)$ and $e_i, \tilde{e}_i \in
\Gamma(E_2)$, we have
\begin{equation}\label{decompbrakcets}
\Phi\circ \{ V, W\}_1 = \sum V^i W^j (\{ e_i, \tilde{e}_j\}_2\circ
\phi) + \sum  {\mathcal L}_{\rho_{*1} V} W^i (\tilde{e}_i
\circ\phi)
 - \sum {\mathcal L}_{\rho_{*1} W} V^i (e_i \circ\phi)
\end{equation}
\end{defi}

It is clear that relations (\ref{decompof}) and
(\ref{decompbrakcets}) are in $\phi^*E_2$. It may appear that
there are ambiguities in this definition. However it can be shown
that the right-hand side of (\ref{decompbrakcets}) is independent
of the $\Phi$-decompositions of $V$ and $W$, for further details
see \cite{higgins}.

\begin{defi}Let  $(E, \M, \rho, \{\,\,,\,\,\})$ be a Lie algebroid.
Then {\it a Lie
  subalgebroid} of $E$ is a morphism of Lie algebroids $\Phi: E' \rightarrow E$,
 $\phi: {\mathcal D} \rightarrow \M$ such that $\Phi$ and $\phi$ are injective
 immersions.
 \end{defi}

In local coordinates we can describe a Lie subalgebroid as
follows. In the neighborhood of a point $x \in {\mathcal D}$ (we
identify $\mathcal D$ with a submanifold of $\M$ ) we choose
coordinates $X^\mu =(X^{\hat{\mu}}, X^{\tilde{\mu}})$ adapted to
the submanifold ${\mathcal D}$ such that in this neighborhood the
submanifold is given by the condition $X^{\hat{\mu}} =0$. We use
the Greek lower case letters with a hat for the coordinates
transverse to the submanifold ${\mathcal D}$ and the same letters
with tilde for the coordinates along the submanifold ${\mathcal
D}$. As well we can introduce the basis on the fiber adapted to
the fact that $E'_x \subset E_x$, namely $e^A = (e^a , e^n)$. We
use the Latin lower case letters from the beginning of alphabet
for the basis of $E'_x$ and from the middle for the remaining
elements in the basis. Then one can show that the above definition
implies the following properties for the anchor map and for
``structure constants'' along ${\mathcal D}$
\begin{equation} \rho^{{\hat{\mu}} a} (0, X^{\tilde{\mu}}) =
0,\,\,\,\,\,\,\,\,\,\,\,\,\,\,\,\,\,
 f^{ab}_{\,\,\,\,\,\,\,\,\,n} (0, X^{\tilde{\mu}}) = 0.
\label{suballocal}
\end{equation}
Thus $\rho^{\tilde{\mu} a} (0, X^{\tilde{\mu}}) $ and
$f^{ab}_{\,\,\,\,\,\,\,\,\,c} (0, X^{\tilde{\mu}})$ define the
structure of a Lie algebroid over $E'\rightarrow{\mathcal D}$.

\subsection{Lie groupoids}

A {\it groupoid} is a small category $\G$ with all arrows
invertible. If the set of objects (points) is $\M$, we say that
$\G$ is a groupoid over $\M$. We shall denote by the same letter
$\G$ the space of arrows, and write
$$
\begin{array}{c}
  \G\\
{\scriptstyle s}\Big\downarrow\Big\downarrow {\scriptstyle t}\\
  \M
\end{array}
$$
where $s$ and $t$ are the source and target maps. If $g,h\in\G$ the
product $gh$ is defined only for pairs $(g,h)$ in the set of
composable arrows
$$  \G^{(2)}=\{(g,h)\in\G\times\G| t(h)=s(g)\},$$
and we denote by $g^{-1}\in \G$ the inverse of $g$, and by ${\rm
id}(x)\equiv x$ the identity arrow at $x\in \M$. The objects $\M$
are thus embedded in $\G$ with $\id$; when no confusion arises we
will omit $\id$ and simply consider $\M\subset\G$. If $\G$ and
$\M$ are topological spaces, all the maps are continuous, and $s$
and $t$ are open surjections, we say that $\G$ is a {\it
topological groupoid}. A {\it Lie groupoid} is a groupoid where
the space of arrows $\G$ and the space of objects $\M$ are smooth
manifolds, the source and target maps $s, t$ are submersions, and
all the other structure maps are smooth. We require $\M$ and the
$s$-fibers $\G_x= s^{-1}(x)$, where $x\in\M$, to be Hausdorff
manifolds, but it is important to allow the total space $\G$ of
arrows to be non-Hausdorf.

The action of $\G$ on a space $X$ equipped of an anchor $\mu: X
\rightarrow \M$ consists in a map from $\G*X=\{(g,x)\in\G\times
X\,|\, s(g)=\mu(x)\}$ to $X$, $(g,x)\rightarrow gx$ such that:
$i)\ \mu(gx)=t(g)$, $ ii)\ g(hx)=(gh)x$, $ iii)\ \mu(x)x = x$.

Given a Lie groupoid $\G$ we can define its tangent Lie algebroid
$\A(\G)$ as follows. It is defined on $\A(\G)_x=T_x\G_x$, for
$x\in \M$; the anchor is $\rho=dt:T_x\G_x\rightarrow T_x\M$. The
bracket comes from the identification of $\Gamma(\A(\G))$ with
left invariant vector fields on $\G$ by choosing the bracket of
vector fields on $\G$. Since not any Lie algebroid come out in
this way, we say that an algebroid $E$ is {\it integrable} if
there exists a Lie groupoid $\G$ such that $\A(\G)=E$. The problem
of integration of Lie algebroid is a generalization of the problem
of integration of Lie algebras.

A groupoid morphism from a groupoid $\G_1$ to $\G_2$ is a
covariant functor; more explicitly, we get the following
definition:

\begin{defi}
\label{def_gpd_morph} Let $\G_i$, $i=1,2$, be Lie groupoids and
let $\M_i$,$\id_i$,$s_i,t_i$ be the corresponding space of units,
their embedding, the source and target maps. A groupoid morphism
from $\G_1$ to $\G_2$ is a couple of maps $(X,\hat{X})$,
$X:\M_1\rightarrow\M_2$ and $\hat{X}:\G_1\rightarrow\G_2$ such
that
\begin{itemize}
\item[{\it i})]
$X\circ s_1=s_2\circ\hat{X}$, $X\circ t_1=t_2\circ\hat{X}$;
\item[{\it ii})] $\hat{X}(ab)=\hat{X}(a)\hat{X}(b)$ for all $a,b\in\G_1^{(2)}$;
\item[{\it iii})] $\hat{X}\circ\id_1=\id_2\circ X$.
\end{itemize}
\end{defi}

Let $E_i=\bigcup\limits_{x\in\M_i} T_{\id_i(x)}(\G_i)_x$ be the tangent
Lie algebroids. It is a fundamental fact that $(X,j)$, where
$j=\hat{X}_*:(E_1)_x\rightarrow(E_2)_{X(x)}$ is a Lie algebroid
morphism.

We have the following Lie theorems for algebroids.

\begin{Thm}\label{LIE1} [{\it Lie I}] Let $E=\A(\G)$ and let $\Ep\subset E$ be a
subalgebroid. Then there exists a Lie subgroupoid $\G'\subset\G$
such that $E'=\A(\G')$.
\end{Thm}

\medskip
\begin{Thm}\label{LIE2} [{\it Lie II}]
Let $\G_1$ and $\G_2$ be two Lie groupoids with Lie algebroids
$E_1$ and $E_2$; if $\G_1$ is source simply connected (ssc) then
for every Lie algebroid morphism $(X,j): E_1\rightarrow E_2$ there
exists a unique groupoid morphism $(X,\hat{X}):\G_1\rightarrow
\G_2$, such that $j=\hat{X}_*$.
\end{Thm}

Moreover for a given
integrable Lie algebroid there exists a unique source simply
connected Lie groupoid integrating it. Next we give
 some basic examples of Lie groupoids.

\begin{Exa}\label{example_group1}
A finite dimensional Lie algebra $\g$, considered as a Lie
algebroid over a point as in Example \ref{example_algbd1}, is
integrated by the simply connected group $G$ seen as a groupoid
over a point.
 \end{Exa}
\begin{Exa}\label{example_group2} The source simply connected groupoid integrating $T \M$
 (see Example \ref{example_algbd2}) is the
fundamental groupoid $\Pi(\M)$, the set of curves in $\M$ modulo
homotopies with fixed end points; the groupoid structures are the
obvious ones, e.g. source (resp. target) is the initial (resp.
final) point, multiplication is concatenation, and identities are
the trivial loops. For each $m\in\M$, $\Pi(\M)_m$ is diffeomorphic
to the universal cover $\tilde \M$ of $\M$ and $\Pi(\M)_m^m$ is
$\pi_1(\M,m)$.
 \end{Exa}
\begin{Exa}\label{example_group3}If the action of $\g$ on $\M$ of Example \ref{example_algbd3}
comes from an action of $G$, where $\g={\rm Lie} G$, then the
action Lie algebroid on $\M\times\g$ is integrated by the action
Lie groupoid $\M\times G$.
 \end{Exa}
\begin{Exa}\label{example_group4} If the algebroid $T^*\M$ associated to a Poisson manifold
 (see Example \ref{example_algbd4})  is
integrable, then the groupoid is a symplectic manifold, called the
{\it symplectic groupoid}. The particular case of a Poisson Lie
group is always integrable. In the factorizable case, the groupoid
integrating it corresponds to the action groupoid $G\times G^*$,
where $G^*$ is the dual Poisson-Lie group, acting on $G$ with the
dressing transformations, see \cite{LW}.  \end{Exa}

We close this section by defining the admissible sections of a
groupoid. For any Lie groupoid $\G$ the group of {\it admissible
sections} $\bis(\G)$ is the group of bisections $\sigma:
\M\rightarrow\G$, such that $s\sigma=\id$ and
$t\sigma=\psi_\sigma:\M\rightarrow \M$ is a diffeomorphism. The
group law is
$\sigma_1\sigma_2(x)=\sigma_1(t(\sigma_2(x)))\sigma_2(x)$ and the
identity is $\id$. Since we have that
$\psi_{\sigma_1\sigma_2}=\psi_{\sigma_1}\psi_{\sigma_2}$,
$\psi_\sigma$ defines an action of $\bis(\G)$ on $\M$ that
preserves the leaves. It comes out that $\bis(\G(E))$ is a Lie
group whose Lie algebra is $\Gamma(E)$. Indeed a tangent vector to
$\sigma\in\bis(\G(E))$ means to assign an element of
$T_{\sigma(x)}\G_x$ for each $x\in \M$; in particular the tangent
space to the identity $\sigma=\id$ is $\Gamma(E)$.

\section{Lie algebroid and TFT}\label{algebroids_TFT}

To any Lie algebroid $(E, \M, \rho, \{\,\,,\,\,\})$ we can
associate a gauge theory in the following way. Consider the space
of bundle maps from the tangent bundle $T\Sigma$ of a
two-dimensional oriented manifold $\Sigma$, possibly with
boundary, to the vector bundle $E$ with base manifold $\M$. We
describe such bundle map by a pair $(X, j)$,
\begin{equation}
\begin{CD}
 T\Sigma @>j  >>   E     \\
@V    VV   @V  VV \\
 \Sigma @>X >> \M
\end{CD}
\end{equation}
where $X : \Sigma \rightarrow \M$ is the base map and $j$ is the
map between fibers, e.g. $j$ is a section in
$\Gamma(T^*\Sigma\otimes X^* E)= \Omega^1(\Sigma, X^* E)$. Now we
consider on $T\Sigma$ the tangent algebroid and we require that
$(X,j)$ is a Lie algebroid morphism. In local coordinates
$\{X^\mu\}$ on $\M$ and $\{u^\alpha\}$ on $\Sigma$ and by choosing
a local trivialization $e^A$ for $E$, $X$ is given by ($\dim \M$)
functions $X^\mu(u)$ and $j$ by ($\rank E$) differential 1-forms
$j_A = j_{A\alpha} du^\alpha$. For arbitrary vector fields on
$\Sigma$, $V,W \in \Gamma(T\Sigma)$ we get the $j$-decomposition
\begin{equation}\label{decompsjuv}
 j\circ V = j_{A\alpha} V^\alpha (e^A\circ X),\,\,\,\,\,\,\,\,\,\,\,\,\,\,\,\,\,
 j\circ W = j_{A\alpha} W^\alpha (e^A\circ X) \;.
\end{equation}
Applying the definition \ref{def_algbd_morph}, we can write
(\ref{decompbrakcets}) in local coordinates
$$ j_{A\alpha} [V, W]^\alpha (e^A\circ X) = j_{A\alpha} V^\alpha j_{B\beta} W^\beta
 (\{e^A, e^B\}\circ X)
+V^\alpha \d_\alpha (j_{A\beta} W^\beta) (e^A\circ X) - $$
\begin{equation}\label{derivmorph}
- W^\alpha \d_\alpha (j_{A\beta} V^\beta) (e^A\circ X)
\end{equation}
where $[\,\,,\,\,]$ is the standard Lie bracket on
$\Gamma(T\Sigma)$ and $\{\,\,,\,\,\}$ is the bracket on
$\Gamma(E)$. The equation (\ref{derivmorph}) implies
\begin{equation}\label{finaleqmorph}
 V^\alpha W^\beta (\d_\beta j_{A\alpha} - \d_\alpha j_{A\beta} -   f^{BC}_{\,\,\,\,\,\,\,\,\,A}
 j_{B\alpha} j_{C\beta}) (e^A\circ X) =0
\end{equation}
where we have used (\ref{definacst}). To summarize the equations
of motion of the system are given by
\ber \label{algeq1} dj_A +
\frac{1}{2} f^{BC}_{\,\,\,\,\,\,\,\,\,A}(X) j_B \wedge j_C = 0,\\
\label{algeq2} dX^\mu - \rho^{\mu A}(X) j_A =0,
\eer{eqsal}
 where
the last equation is a simple consequence of (\ref{commancmap}).

 Thus  the fact that $(X, j)$ is a Lie algebroid morphism implies the first order differential
 equations for $X^\mu$ and $j_{A\alpha}d\xi^\alpha$.
 These equations form a consistent system of partial differential equations due to the
  properties (\ref{consalg1}) and (\ref{consalg2}). On the top of this the system is invariant under
 the infinitesimal gauge transformations
\ber
\label{gaualg1} \delta j_A= -d\beta^A -f^{BC}_{\,\,\,\,\,\,\,\,\,A}(X)j_B \beta_C \\
\label{gaualg2} \delta X^\mu = - \rho^{\mu A} (X)\beta_A
\eer{gau}
 where $\beta$ is a gauge parameter. For the similar
 discussion of the system  (\ref{algeq1})-(\ref{gaualg2})  see \cite{Bojowald:2004wu}.

 The motivating example for the system (\ref{algeq1})-(\ref{gaualg2}) is PSM where
 $\Sigma$ is two-dimensional and $E= T^* \M$ defined in Example \ref{example_algbd4}. In this case the equations
   (\ref{algeq1})-(\ref{algeq2}) are the stationary points of the action functional
   \begin{equation}\label{PSMaction}
   S(X,\eta) = \int\limits_\Sigma  \langle \eta,  dX\rangle +
 \frac{1}{2} \langle \eta, (\alpha \circ X) \eta \rangle ,
    \end{equation}
    where $(X, \eta)$ is bundle morphism from $T\Sigma$ to $T^*\M$.  The pairing $\langle , \rangle$
     is defined as pairing of the values in $T\M$ and $T^*\M$
       and the exterior product of differential forms. The action functional is invariant under
        the corresponding transformations (\ref{gaualg1})-(\ref{gaualg2}).

  Indeed the system  (\ref{algeq1})-(\ref{gaualg2})
  is defined for $\Sigma$ being a manifolds of any dimension, $\dim \Sigma =p$.
  The equations of motion can be derived from the following trivial action
\begin{equation}\label{lagranmultact}
 S (X, j, A, \lambda)
  = \int\limits_{\Sigma} \,\lambda^A\wedge ( d j_A + \frac{1}{2} f^{BC}_{\,\,\,\,\,\,\,\,\,A}(X) j_B \wedge
 j_C) + A_\mu \wedge (dX^\mu - \rho^{\mu A} (X) j_A),
\end{equation}
 where $\lambda \in \Omega^{p-2}(\Sigma, X^* E^*)$
  and $A  \in \Omega^{p-1}(\Sigma, X^* T^* \M)$ are the Lagrangian multipliers, see \cite{strobl}
  for an analogous discussion.
 The action (\ref{lagranmultact}) is invariant under
(\ref{gaualg1}) and (\ref{gaualg2}) together with the additional
transformations
\begin{eqnarray}\label{transforB}
 \label{GGGG1}&&\delta \lambda^A = \rho ^{\mu A} b_\mu +f^{AB}_{\,\,\,\,\,\,\,\,\,C}\lambda^C \beta_B \\
\label{GGGG2} && \delta A_\mu =(-1)^{(p -1)}d b_\mu -(\d_\mu f^{AB}_{\,\,\,\,\,\,\,\,\,C} )\lambda^C \wedge j_A \beta_A
 - (\d_\mu \rho^{\nu A})   b_\nu \wedge j_A + (\d_\mu \rho^{\nu A}) A_\nu \beta_A ,
\end{eqnarray}
 where $b\in\Omega^{p-2}(\Sigma,X^*T^*M)$ and $\beta\in\Omega^0(X^*E)$ are the gauge parameters.
 In coordinate free way the fields
  $(X, j, A, \lambda)$ can be interpreted as follows. $(X, j)$ is a bundle morphism $T\Sigma \rightarrow
   E$, $(X, A)$ is a bundle morphism $\wedge^{p-1} T\Sigma \rightarrow T^* \M$ and
    $(X, \lambda)$ is a bundle morphism $\wedge^{p-2} T\Sigma \rightarrow E^*$.  Introducing
     the pairings as pairing in values between $E$ and $E^*$ and as paring between $T\M$ and
      $T^*\M$ and the exterior product of differential forms on $\Sigma$ we can rewrite the action
       functional (\ref{lagranmultact}) in coordinate independent form.  This theory is a rather obvious
        generalization of $BF$-theory \cite{Horowitz:1989ng, Blau:1989bq} to the case of a generic Lie algebroid.

 If $\Sigma$ is two dimensional then the action
 (\ref{lagranmultact}) can be written in the following form
\begin{equation}\label{bigPoisson}
 S= \int\limits_{\Sigma} j_A \wedge d\lambda^A + A_\mu \wedge dX^\mu + \frac{1}{2} ( f^{BC}_{\,\,\,\,\,\,\,\,\,A}(X)
 \lambda^A) j_B \wedge j_C + \rho^{\mu A}(X) j_A \wedge A_\mu
\end{equation}
 and it differs from (\ref{lagranmultact}) only by a boundary term. The action (\ref{bigPoisson}) has
 a clear interpretation, this is the Poisson sigma model for $E^*$. In fact, if
 $(E, \M, \rho, \{\,\,,\,\,\}) $ is a Lie algebroid, the dual bundle $E^*$ has a natural
 Poisson structure, defined by the tensor $\pi\in\wedge^2T^*E^*$
 given in coordinates $X^\mu,\lambda^A$ as
\begin{equation}\label{defpoisEstar}
 \pi (X, \lambda)  = f^{AB}_{\,\,\,\,\,\,\,\,\,C} \lambda^C \frac{\d}{\d\lambda^A}\wedge
 \frac{\d}{\d \lambda^B} + \rho^{\mu A} \frac{\d}{\d \lambda^A} \wedge \frac{\d}{\d X^\mu} .
\end{equation}
In the action (\ref{bigPoisson}) we have that $(X,\lambda) :
\Sigma \rightarrow E^*$ and $(A, j)$ is a differential form on
$\Sigma$ taking values in the pull-back by $(X,\lambda)$ of
$T^*E^*$.

Every solution of (\ref{algeq1}) and (\ref{algeq2}) defines a
solution of the equations of motion of (\ref{bigPoisson}), e.g. we
have embedded our system of equations in a TFT. Let us see it in
an intrinsic way. The full set of equations, including those
obtained by varying $X$ and $j$, obviously describes the algebroid
morphisms from $T\Sigma$ to $T^*E^*$.  In Lemma 4.2 of
\cite{alberto1} it is shown that $E$ is a subalgebroid of
$T^*E^*$. The injection is defined as follows: the fibre
$T^*_{(m,\alpha)}E^*$ over $(m,\alpha)\in\M\times E_m$, is
$T^*_m\M\oplus E_m$. The injection is defined as
$\iota:E_m\rightarrow T^*_{(m,0)}E^*$ as
$\iota(m,a)=((m,0),0\oplus a)$. It comes out that this is an
injective algebroid morphism. Thus composing with $\iota$ we
can inject the set of vector bundle morphisms from $T\Sigma$ to
$E$ into the space of fields of the model and every  algebroid
morphism of $E$ defines an algebroid morphism for $T^*E^*$, e.g.
is a solution of the equations of motion of the PSM with target
$E^*$. Of course this mapping is not surjective, e.g. there are
solutions with $\lambda$ and $A$ different from $0$.

If $\d \Sigma \neq \emptyset$ we have to choose  appropriate
boundary conditions on the fields. For example, we may ask about
that the boundary terms vanish in
  the variations of (\ref{bigPoisson}).  Thus
 in order to get the equations of
motion, we have to choose boundary conditions such that
\begin{equation}\label{boundcgenv}
 (j_{\tau A} \delta \lambda^A + A_{\tau \mu} \delta X^\mu) |_{\d \Sigma} = 0
\end{equation}
 where $j_A|_{T^* \d\Sigma} = j_{A\tau} d\tau$ and $A_\mu|_{T^*\d \Sigma} =
 A_{\mu\tau} d\tau$. The action (\ref{bigPoisson})  is invariant under the infinitesimal gauge
transformations (\ref{gaualg1})-(\ref{gaualg2}) and (\ref{GGGG1})-(\ref{GGGG2}) with
parameter $(b_\mu dX^\mu , \beta_A d\lambda^A) \in \Gamma
 ((X,\lambda)^*(T^* E^*))$, provided the following boundary condition is
 satisfied
\begin{equation}\label{boundcongaug}
 (\beta_A \d_\tau \lambda^A + b_\mu \d_\tau X^\mu)|_{\d \Sigma} =
 0\;.
\end{equation}
Finally, also the boundary conditions should be invariant under
the residual gauge transformations. By using the results of
\cite{CF2} (see also \cite{bonechi}) one can establish that the
boundary conditions for the theory are labeled by the coisotropic
submanifolds of $E^*$. Recall that a submanifold $X$ of a Poisson
manifold $E^*$ is coisotropic iff the conormal bundle $N^*X$ is a
subalgebroid of $T^*E^*$.

Motivated by this discussion it is natural to choose those
boundary conditions that come from $E$. In fact for every
subalgebroid $E'\subset E$ over $\D\subset\M$, we can see that
$E'{}^\perp\subset E^*$ is a coisotropic submanifold. Let us
choose the adapted coordinates $X^\mu = (X^{\hat{\mu}},
X^{\tilde{\mu}})$ and trivialization  $e_A= (e_a, e_n)$ such that
the $\D$ corresponds to $X^{\hat{\mu}}=0$ and $E'_X=\langle
e_a\rangle$. It is then easy to verify that the local conditions
(\ref{suballocal}) for $E'$ be a subalgebroid correspond to those
for $E'{}^\perp$ be a coisotropic submanifold, $\pi^{ab} =0$ and
$\pi^{a\hat{\mu}}=0$. The contrary is not true as can be easily
understood from the following example.

\begin{Exa} Let $E=\g$ be a Lie algebra and $\g^*$ is vector space with the canonical Poisson
 structure. Then the subalgebroids of $\g$ are the Lie
subalgebras $\h\subset\g$ that define the linear coisotropic
submanifolds $\h^\perp\subset\g^*$. However not any coisotropic
submanifold of $\g^*$ arises in this way.
\end{Exa}

Thus in forthcoming discussion when we refer to the open case we consider the
 system (\ref{algeq1})-(\ref{gaualg2})  with the boundary conditions
given by  Lie subalgebroids of  $(E, \M, \rho,
\{\,\,,\,\,\})$
\begin{equation}
\begin{CD}
 T\d\Sigma @>j  >>   E' @> \Phi >> E    \\
@V    VV   @V  VV  @V VV \\
 \d\Sigma @>X >> {\mathcal D} @>\phi >> \M
\end{CD}
\end{equation}
 On the boundary the gauge transformations are restricted correspondently.

 We close this Section with a comment about PSM with target $\M$ when $\M$ is a Poisson manifold.
 The boundary conditions for this PSM are not defined by any subalgebroid of $T^*\M$, but only by those which
are conormal bundles of a submanifold  \cite{CF2}. This is not
surprising since we have motivated our choice of boundary
conditions starting from PSM with target $(T^*\M)^*$. It is not
clear at the moment the relevance of this wider class of boundary
conditions in the context of the PSM with target $\M$.

\bigskip
\bigskip

\section{Integration of Gauge transformations}
\label{GAUGE}

In this section we define the finite gauge transformations that
integrate the infinitesimal transformations (\ref{gaualg1}) and
(\ref{gaualg2}). In the case of an integrable Lie algebroid, we
will analyze groupoid morphisms rather than algebroid morphisms,
since it is much easier to introduce the finite gauge
transformations.

In fact, due to Theorem \ref{LIE2} every solution $(X,j)$ of the
equations (\ref{algeq1}), (\ref{algeq2}) can be lifted to a
groupoid morphism $(X,\hat{X})$ between $\Pi(\Sigma)$ and $\G(E)$,
the (ssc) groupoid integrating $E$, and vice-versa. In the
following we will identify the solutions of the equations
(\ref{algeq1}) and (\ref{algeq2}) with the groupoid morphisms
from $\Pi(\Sigma)$ to $\G(E)$ that they generate. We denote
the space of all morphisms from $\Pi(\Sigma)$
 to $\G(E)$ with ${\rm Mor}(\Pi(\Sigma),\G(E))$.

\subsection{The closed case}
\label{gaugeclosed}

In this subsection we consider the case when $\d \Sigma =\emptyset$.
We assume that $E$ is integrated by $\G(E)$ with
$s,t:\G(E)\rightarrow \M$ being the source and target map. Let
$\id:\M\rightarrow\G$ be the usual embedding of $\M$ in $\G(E)$ as
the space of identities. As usual we denote with $\G(E)_x$
($\G(E)^x$) the fiber of the source (target) map in $x\in\M$ and
with $\G(E)_x^y=\G(E)_x\cap\G(E)^y$. Recall that $\G(E)_x$ and
$\G(E)^x$ are separable smooth manifolds and that
$T_x\G(E)_x=E_x$.

We will first define the gauge transformations on the morphisms
 from $\Pi(\Sigma)$ to $\G(E)$ and then we will compute the
induced transformations on the algebroid morphisms between
$T\Sigma$ and $E$.

Following \cite{BCZ} we introduce the infinite-dimensional
groupoid $\G^\Sigma = \{\hat{\Phi}\colon \Sigma\rightarrow\G(E)\}$
over $\M^\Sigma=\{\Phi\colon \Sigma\rightarrow \M\}$ with
structure maps defined pointwise. Namely, we define source and
target by $s(\hat{\Phi})(u)=s({\hat\Phi}(u))$,
$t(\hat{\Phi})(u)=t({\hat\Phi}(u))$ for $u \in \Sigma$ and
multiplication by
 $\hat{\Phi}_1\hat{\Phi}_2 (u) =
\hat\Phi_1(u)\hat\Phi_2(u)$. A section $S$ of the associated
algebroid\footnote{This algebroid has been defined
for one dimensional $\Sigma$  in \cite{BC}. It has been done intrinsically in
terms of the Lie algebroid $E$ and thus it exists also for
nonintegrable algebroids.} $A(\G^\Sigma)$  is defined by giving a section
$S(\Phi) \in \Gamma(\Phi^*E)$ for every $\Phi\in \M^\Sigma$.
There is a natural groupoid action of $\G^\Sigma$ on ${\rm
Mor}(\Pi(\Sigma),\G(E))$ which is given by
\begin{equation}
X_{\hat\Phi}(u)=t(\hat\Phi)(u) \;\;\;\; \X_{\hat\Phi}([c_{uv}])
=\hat\Phi(u)\X([c_{uv}])\hat\Phi(v)^{-1}, \label{gauge_trans_gpd}
\end{equation}
where $(X, \X), (X_{\hat\Phi}, \X_{\hat\Phi})\in{\rm
Mor}(\Pi(\Sigma),\G(E))$, $\hat\Phi\in\G^\Sigma$ with
$s(\hat\Phi)=X$ and $[c_{uv}]$ is the homotopy class of a curve
$c_{uv}$ in $\Sigma$.  Thus we declare $\G^\Sigma$ as our choice
of finite gauge transformations.

However there are alternative choices of finite gauge
transformations, e.g. the group $\bis (\G^\Sigma)$ of bisections
of $\G^\Sigma$. In this case the  formula (\ref{gauge_trans_gpd})
also defines a group action of $\bis(\G^\Sigma)$ on ${\rm
Mor}(\Pi(\Sigma),\G(E))$. The orbits of $\G^\Sigma$ contain the
orbits of $\bis(\G^\Sigma)$ and thus the choice of $\G^\Sigma$ is
more generic one.  Another possible choice of the gauge
transformations is $(\bis\G(E))^\Sigma$, which is the group of
maps from $\Sigma$ to $\bis \G(E)$, seen as a subgroup of
$\bis(\G^\Sigma)$, see section 3.1 in \cite{BCZ}. However it is
very hard to work with these groups and in the following we will
consider only the groupoid action of $\G^\Sigma$ on ${\rm
Mor}(\Pi(\Sigma),\G(E))$.

\bigskip

Indeed the choice of $\G^\Sigma$ as finite gauge transformations
looks natural from the categorical point of view. Namely if we
look to the above groupoids as categories then any groupoid
morphism in ${\rm Mor}(\Pi(\Sigma),\G(E))$ is a covariant functor
from $\Pi(\Sigma)$ to $\G(E)$ and a gauge transformation between
two groupoid morphisms as defined in (\ref{gauge_trans_gpd}) is a
natural transformation between the functors. In Section \ref{END}
we will comment more on this issue.

\bigskip

 There is another possibility to introduce the notion of gauge equivalence
  between the groupoid morphisms or algebroid morphisms, e.g. see \cite{Bojowald:2004wu}.
   Namely this can be done via groupoid (algebroid) homotopies.
 The  groupoid (algebroid) homotopy is an
alternative way of integrating the gauge transformations
(\ref{gaualg1}) and (\ref{gaualg2}). Let $I=[0,1]$; it is clear
that for groupoid $\Pi(\Sigma)$ over $\Sigma$, we can define on $\Pi(\Sigma)\times
I\times I$ a groupoid structure over $\Sigma \times I$ with the corresponding algebroid
 given by $T(\Sigma \times I)$.

\smallskip
\begin{defi} \label{gpd_hom}
Let $\hat{X}_i$, $i=1,2$, be two groupoid morphisms from
$\Pi(\Sigma)$ to $\G(E)$. We say that $\X_1$ and $\X_2$ are {\it
homotopic} if there exists a groupoid morphism
$\hat{X}_{12}:\Pi(\Sigma)\times I\times I \rightarrow\G(E)$ such
that $\hat{X}_{12}(-,0,0)=\hat{X}_1$ and
$\hat{X}_{12}(-,1,1)=\hat{X}_2$.
\end{defi}
\smallskip
\begin{defi} \label{alg_hom}
Let $(X_i, j_i)$, $i=1,2$, be two algeboid morphisms from
$T\Sigma$ to $E$. We say that $(X_1, j_1)$ and $(X_2, j_2)$ are
{\it homotopic} if there exists an algebroid morphism $ (X_{12},
j_{12}) : T(\Sigma \times I)
 \rightarrow E$ such that $(X_{12}, j_{12})(-,0)=(X_1, j_1)$ and
$(X_{12}, j_{12})(-,1)=(X_2, j_2)$.
\end{defi}

 Next we show that the groupoid homotopies are the gauge
transformations connected to the identities, {\it i.e.} those
transformations that live in the component $(\G^\Sigma_X)_o$ of
the source fibre $\G^\Sigma_X$ over $X:\Sigma\rightarrow \M$
connected to the identity $X$. Borrowing the terminology from
gauge theory we say that groupoid homotopies are the small gauge
transformations of $\G^\Sigma$.

\begin{Lem}\label{homotopies_gauge}
Two groupoid morphisms  $\hat{X}_{i}:\Pi(\Sigma)\rightarrow\G(E)$,
$i=1,2$, are homotopic if and only if there exists a gauge
transformation $\hat\Phi\in(\G^\Sigma_{X_1})_o$ such that
$\hat{X}_2=(\hat{X}_1)_{\hat\Phi}$.
\end{Lem}

{\it Proof}. Let $\hat{X}_i$ be homotopic with homotopy
$\hat{X}_{12}$. We have that
\begin{eqnarray*}
\X_2[c_{uv}]&=&\X_{12}([c_{uv}],1,1)=
\X_{12}([c_{uu}^{tr}],1,0)\X_{12}([c_{uv}],0,0)\X_{12}([c_{vv}^{tr}],0,1)\\
&=& \hat\Phi(u) \X_1([c_{uv}]) \hat\Phi(v)^{-1} \;,
\end{eqnarray*}
where $[c_{uu}^{tr}]$ is the class of the trivial loop through
$u$, $\hat\Phi\in\G^\Sigma$ is defined by $\hat\Phi(u)=
\X_{12}([c_{uu}^{tr}],1,0)$ and
$\hat{\gamma}(s)(u)=\hat{X}_{12}([c^{tr}_{uu}],s,0)\in\G(E)_{X_1(u)}$
is such that $\hat{\gamma}(0)=X_1$ and
$\hat{\gamma}(1)=\hat{\Phi}$.

Conversely, let $\hat{\Phi}\in(\G^\Sigma_{X_1})_o$, with
$\hat{\gamma}:I\rightarrow\G^\Sigma_{X_1}$ such that
$\hat{\gamma}(0)=X_1$ and $\hat{\gamma}(1)=\hat{\Phi}$, and
$\hat{X}_2=(\hat{X}_1)_\Phi$. The required homotopy is then
defined as
$\hat{X}_{12}([c_{uv}],s_1,s_2)=\hat{\gamma}(s_1)(u)\hat{X}_1[c_{uv}]\hat{\gamma}(s_2)(v)^{-1}$.
\endproof

Below for completeness we describe the action of $\G^\Sigma$ on
algebroid morphisms.

\medskip
\subsection{Groupoid action on the algebroid morphisms}
\label{groupact}

In this subsection we compute the groupoid action of $\G^\Sigma$
on the algebroid morphism $(X,j):T\Sigma\rightarrow E$. Let
$\hat\Phi\in\G^\Sigma$ be such that $s(\hat\Phi)=X$ and let
$\X:\Pi(\Sigma)\rightarrow \G(E)$ be the groupoid morphism
integrating $(X,j)$, e.g.$j=\X_*$. We define the action of
$\hat\Phi$ on $j$ as $j_{\hat\Phi}=\X_{\hat\Phi}{}_*$.

Let $\G(E)^{(2)}_x=\{(\gamma_1,\gamma_2)\in\G(E)\times\G(E)_x\,
|\, s(\gamma_1)=t(\gamma_2)\}$ and let $m:\G(E)^{(2)}_x\rightarrow
\G(E)_x$ be the multiplication
$m(\gamma_1,\gamma_2)=\gamma_1\gamma_2$. For each $v\in\Sigma$ we
have that
$\X_{\hat\Phi}:\Pi(\Sigma)_v\rightarrow\G(E)_{t(\hat\Phi(v))}$ can
be expressed as $ \X_{\hat\Phi}= R_{\hat\Phi(v)^{-1}}\circ m \circ
(\hat\Phi \circ t,\X)$, where $(\hat\Phi\circ
t,\X):\Pi(\Sigma)_v\rightarrow\G(E)^{(2)}_{X(v)}$, and
$R_{\hat\Phi(v)^{-1}}:\G(E)_{X(v)}\rightarrow\G(E)_{t(\hat\Phi(v))}$
denotes the right multiplication by $\hat\Phi(v)^{-1}$. The
tangent map $j_{\hat\Phi}:T_v\Sigma\rightarrow E_{t(\hat\Phi(v))}$
is expressed on $w\in T_v\Sigma$ as
\begin{equation}\label{grpd_action_algbd_morph}
j_{\hat\Phi}(w) = R_{\hat\Phi(v)^{-1}}{}_*\circ
m_*(\hat\Phi_*(w)\oplus j(w)).
\end{equation}

Remark that in (\ref{grpd_action_algbd_morph}) we consider
$\hat\Phi_*(w)\oplus j(w)\in
T_{(\hat\Phi(v),X(v))}\G(E)^{(2)}_{X(v)}={\rm Ker}(s_*-t_*)\subset
T_{\hat\Phi(v)}\G(E)\oplus T_{X(v)}\G(E)_{X(v)}$, where
$s_*-t_*:T_{\gamma_1}\G(E)\oplus T_{\gamma_2}\G(E)_x\rightarrow
T_{s(\gamma_1)}M$ is defined on every $(\gamma_1,\gamma_2)\in
\G(E)^{(2)}_{s(\gamma_1)}$. It is then clear that
(\ref{grpd_action_algbd_morph}) makes sense only if
$s_*\circ\hat\Phi_*(w)=t_*\circ j(w)$, e.g. $X_*=\rho\circ j$
which is (\ref{algeq2}). So the action of $\G^\Sigma$
automatically extends only to those vector bundle morphisms that
commute with the anchor maps. The correct definition of the off
shell gauge transformations is delicate and is beyond the scope of the present work, e.g.
  see
\cite{Bojowald:2004wu} for a discussion of this problem.

Indeed the above construction is a direct generalization of the
following example.

\begin{Exa}
\label{FLAT11}
 Consider the Examples \ref{example_algbd1} and \ref{example_group1}.
Let $E$ be the Lie algebra $\g$; then $\G(E)=\G(E)^{(2)}=G$ and
$m_*=R_{g_2}{}_*+L_{g_1}{}_*:T_{g_1}\oplus T_{g_2}\rightarrow
T_{g_1g_2}$. If we plug it in (\ref{grpd_action_algbd_morph}) we
get the action of $g\in G^\Sigma$ on $j:T_v\Sigma\rightarrow\g$ as
\begin{equation}\label{flatconnection}
j_g(w)=(R_{g(v)^{-1}}{}_*\circ g_*+\Ad_{g(v)}{}_*)\circ j(w).
\end{equation}
 In this example a Lie algebroid morphism is a flat connection on the trivial bundle and (\ref{flatconnection})
 is the gauge transformation of a connection. The associated groupoid
morphism is defined by the parallel transport which transforms
with the adjoint, accordingly to (\ref{gauge_trans_gpd}).
\end{Exa}

\bigskip

\subsection{The open case}\label{gauge_open}

Let $\Sigma$ be a surface with $n$ boundary components,
$\partial\Sigma=\bigcup\limits_{i=1}^n \partial_i\Sigma$.
According to the discussion in Section \ref{algebroids_TFT}, we
consider a set $\E=\{E_i\}$ of $n$ subalgebroids $E_i\subset E$
over $\D_i$. Due to the integrability of $E$, there are $n$ source
connected Lie subgroupoids $\G(E_i)\subset\G(E)$ that integrate
the Lie subalgebroids $E_i$.

We consider the space ${\rm Mor}(\Pi(\Sigma),\G(E); \G(\E))$ of
groupoid morphisms $\X:\Pi(\Sigma)\rightarrow\G(E)$ such that
${\X}\left(\Pi(\Sigma)|_{\partial_i\Sigma}\right)\subset\G(E_i)$,
where
$\Pi(\Sigma)|_{\partial_i\Sigma}=\{[c_{uv}]\in\Pi(\Sigma)\,|\,
c_{uv}\subset\partial_i\Sigma\}$ is the subgroupoid of
$\Pi(\Sigma)$ integrating the $i$-th boundary component
$\partial_i\Sigma$.

In analogy with the closed case, we define the groupoid
$\G^{\Sigma,\E}$ over $\M^{\Sigma,\E}$, where
$$\G^{\Sigma,\E}=\{\hat{\Phi}\colon \Sigma\rightarrow\G(E)\,|\,
\hat{\Phi}(\partial_i\Sigma)\subset\G(E_i)\}\,,~~\M^{\Sigma,\E}=\{\Phi\colon
\Sigma\rightarrow \M\,|\, \Phi(\partial_i\Sigma)\subset \D_i \}.$$
Formula (\ref{gauge_trans_gpd}) gives a groupoid action of
$\G^{\Sigma,\E}$ on ${\rm Mor}(\Pi(\Sigma),\G(\M);\G(\E))$.

From subsections \ref{gaugeclosed} and \ref{groupact} we can generalize all results for the case
 $\d \Sigma \neq \emptyset$.  The generalizations are rather straightforward.
 Thus we can introduce groupoid homotopies respecting the boundary
conditions.

\begin{defi}
The groupoid morphisms $\hat{X}_i\in{\rm
Mor}(\Pi(\Sigma),\G(E);\G(\E))$, $i=1,2$, are homotopic if there
exists a groupoid morphism $\hat{X}_{12}:\Pi(\Sigma)\times I\times
I \rightarrow\G(E)$ such that $\hat{X}_{12}(_,0,0)=\hat{X}_1$,
$\hat{X}_{12}(-,1,1)=\hat{X}_{12}$ and
$\hat{X}_{12}:\Pi(\Sigma)|_{\partial_i\Sigma}\times I\times I
\rightarrow\G(E_i)$.
\end{defi}

In analogy to the closed case, groupoid homotopies are the
 gauge transformations connected to the identities. In fact,
let $(\G^{\Sigma,\E}_X)_o$ be the component connected to
$X:\Sigma\rightarrow\M$ of the source fiber over $X$. By repeating
the same proof as in  Lemma {\ref{homotopies_gauge} and taking care of
boundary conditions we can prove the following result.

\begin{Lem}\label{homotopies_gauge_boundary}
Two groupoid morphisms $\hat{X}_{i}\in{\rm
Mor}(\Pi(\Sigma),\G(E);\G(\E))$, $i=1,2$, are homotopic if and
only if there exists a gauge transformation
$\hat\Phi\in(\G^{\Sigma,\E}_{X_1})_o$ such that
$\hat{X}_2=(\hat{X}_1)_{\hat\Phi}$.
\end{Lem}

\medskip
%Let us choose now $\Sigma=I$ with boundary conditions defined by
%the subalgebroids $\E=\{E_\alpha\}_{\alpha=0,1}$, on the boundary
%components $\partial I=\{0,1\}$.

%\begin{Lem}
%The gauge groupoid $\G(E)^{I,\E}$ is connected to the identities.
%\end{Lem}

%{\it Proof}. Let us recall that the groupoid $\G(E)$ can be
%constructed as the space of algebroid morphisms $TI\rightarrow E$
%divided algebroid homotopies, {\bf [add citation of cattaneo \cite{SG}
%felder and fernandes crainic \cite{FC}]}. So we identify every
%$\gamma\in\G(E)$ as an equivalence class of an algebroid morphism
%$\gamma=[\xi,\eta]$, with $(\xi,\eta):TI\rightarrow E$. Let us
%define $R:\G(E)\times I\rightarrow\G(E)$ as
%$R([\xi,\eta],t)=[\xi_t,\eta_t]$, where $\xi_t(\tau)=\xi(\tau t)$
%and $\eta_t(\tau)=t\eta(\tau t)$. It is easy to see that $R$ is a
%deformation retract from $\G(E)$ to $\M$, such that $i$)
%$R(\gamma,0)=s(\gamma)$; $ii$) $R(\gamma,1)=\gamma$; $iii$)
%$R(\gamma,t)\in\G(E)_{s(\gamma)}$; $iv$) $R(m,t)=m$, for
%$\gamma\in\G(E)$, $m\in\M$, $t\in I$.

%Then for every $\hat\Psi\in\G(E)^{I,\E}$ we have that
%$c:I\rightarrow\G(E)^{I,\E}$ defined by
%$c(t)(\tau)=R(\hat\Psi(\tau),t)$ is a curve in
%$\G(E)^{I,\E}_{s(\hat\Psi)}$ joining $s(\hat\Psi)$ with
%$\hat\Psi$. \endproof

\bigskip
\bigskip

\section{Moduli Space of Solutions}
\label{MODULI}

In this section we discuss the moduli space of solutions modulo
gauge transformations. As explained in the previous section, among
several choices for the finite gauge transformations, we will
choose the largest set, the groupoid $\G^\Sigma$. The main
motivation for this choice comes from the PSM on the disk. In
\cite{BCZ} it has been shown that this is the correct gauge group
for the observables that are relevant for deformation
quantization. Moreover the discussion of the moduli space is
particularly simple and is the straightforward generalization of
that obtained in \cite{bonechi} for Poisson-Lie groups.

As we will see the whole construction is a quite direct
generalization of the moduli space of flat $G$-connections. Let us
analyze this case first. The algebroid morphisms from $T\Sigma$ to
$\g$ correspond to the flat connections in the trivial $G$-bundle
over $\Sigma$, where $G$ is the simply connected Lie group
integrating $\g$. The meaning of the second Lie theorem is that we
can equivalently describe any flat connection for the trivial
bundle by assigning the parallel transport. If $\dim\Sigma\leq 2$ then
there are no other topologically inequivalent $G$-bundles and therefore
the moduli space of algebroid morphisms coincides with the moduli
space of flat connections; in generic dimension it will be the
component corresponding to the trivial bundle. We are going to
show that this description remains valid if we consider a generic
integrable algebroid $E$ and take the (ssc) groupoid $G(E)$
integrating it.

\subsection{The closed case}\label{moduli_closed}

Let $\Sigma$ be a closed surface. We denote by $\M(\Sigma,\G(E))$
the space of groupoid morphisms divided by the action
(\ref{gauge_trans_gpd}) of $\G^\Sigma$, i.e.
\begin{equation}
\label{moduli_space_closed} \M(\Sigma,\G(E)) = {\rm
Mor}(\Pi(\Sigma),\G(E))/\G^\Sigma\;.
\end{equation}

It is clear that, for any $(X,\hat{X})\in {\rm
Mor}(\Pi(\Sigma),\G(E))$, $X(\Sigma)$ is contained in a single
leaf $L$. Since gauge transformations do not change the leaf, we
can decompose the total moduli space (\ref{moduli_space_closed})
in the union of the moduli space $\M(\Sigma,\G(E);L)$ of solutions
corresponding to the leaf $L$.

The following moduli space will be of central interest for us. Let
us fix a leaf $L\subset\M$ and consider for any point $x_0\in L$
the isotropy group $\G_{x_0}^{x_0}=\G(E)^{x_0}_{x_0}$ together
with the moduli space of flat connections
$F(\Sigma,\G_{x_0}^{x_0})={\rm
Hom}(\pi_1(\Sigma),\G_{x_0}^{x_0})/\Ad$, where $\Ad$ denotes the
adjoint action of $\G_{x_0}^{x_0}$. By varying $x_0\in L$, the
isotropy group changes by conjugation so that the moduli spaces
are isomorphic; we denote it $F(\Sigma,\G(E);L)$. We define the
moduli space of generalized flat $\G$-connections on $\Sigma$ as
the union over all symplectic leaves, e.g. $F(\Sigma,\G(E))=\bigcup\limits_L
F(\Sigma,\G(E);L)$. It is important to notice that the moduli space can
be introduced without referring to any choice of $x_0\in L$. In
fact one can verify that
\begin{equation}\label{generalized_flat}
F(\Sigma,\G(E))=\bigcup\limits_{y\in \M}{\rm
Hom}(\pi_1(\Sigma),\G_y^y)/\Ad \;,
\end{equation}
where $\Ad$ means the adjoint action of $\G(E)$, $\phi\rightarrow
\gamma \phi \gamma^{-1}$, for $\phi\in{\rm
Hom}(\pi_1(\Sigma),\G_y^y)$ and $\gamma\in\G_y$. Equivalently
$F(\Sigma,\G(E))={\rm Mor}(\pi_1(\Sigma),\G(E))/\Ad$ is the
moduli space of groupoid morphisms from $\pi_1(\Sigma)$ to
$\G(E)$, where $\pi_1(\Sigma)$ is regarded as a groupoid over a
point.

\begin{rmk}\label{topology}
{\rm There is a very natural topology on $F(\Sigma,\G(E))$. Let
$\Sigma$ be the compact surface of genus $g$, then $\pi_1(\Sigma)$
is generated by $\{a_i,b_i\}_{i=1}^g$ with relation
$\Pi_i[a_i,b_i]=1$, where $[a,b]=aba^{-1}b^{-1}$. Let us define
$\G^2_2=\bigcup\limits_{x\in\M}\times^{2g}\G(E)^x_x\subset\times^{2g}\G(E)$
together with the map $p:\G^2_2\rightarrow\G(E)$,
$$p(z_1,w_1,\ldots,z_g,w_g) = [z_1,w_1]\ldots[z_g,w_g]\;.$$
It is clear that $p^{-1}(\M)\subset\times^{2g}\G(E)$ inherits the
relative topology and that $F(\Sigma,\G(E))$ $=p^{-1}(\M)/{\rm
Ad}$ the quotient topology.}\endproof
\end{rmk}

In Proposition \ref{prop_moduli_space_closed} we show that the two
moduli spaces (\ref{moduli_space_closed}) and
(\ref{generalized_flat}) coincide. Before proving it, we will
introduce the following auxiliary constructions. Let us fix
$u_0\in\Sigma$; we identify $\Pi(\Sigma)_{u_0}^{u_0}$ as
$\pi_1(\Sigma)$ and $\Pi(\Sigma)_{u_0}$ as $\tilde{\Sigma}$, the
universal cover of $\Sigma$. Let us introduce the following
trivialization for the $\pi_1(\Sigma)$ principal bundle
$\tilde{\Sigma}\rightarrow\Sigma$. Let $\{U_\alpha\}_\alpha$ be a
covering of $\Sigma$ with $U_\alpha$ and $U_\alpha\cap U_\beta$
contractible. Let us fix for each $u\in U_\alpha$ a curve
$c^\alpha_{u_0u}$ starting in $u$ and ending in $u_0$ in such a
way that, once $U_\alpha$ is contracted, all such curves are
homotopic. Then let us define on $U_\alpha\cap U_\beta$,
$h_{\alpha\beta}=c^{\alpha}_{u_0u}\circ c^{\beta}_{uu_0}$; it is
clear that $[h_{\alpha\beta}]\in\pi_1(\Sigma)$ is constant for all
$u\in U_\alpha\cap U_\beta$; we have defined a flat structure on
$\tilde{\Sigma}$.

\medskip

\begin{pro}\label{prop_moduli_space_closed}
For every closed manifold $\Sigma$ such that $\dim\Sigma=1,2$ and for
every source simply connected Lie groupoid $\G(E)$ we have
$\M(\Sigma,\G(E))= F(\Sigma,\G(E)) \;.$
\end{pro}

{\it Proof}. For any $(X,\X)\in{\rm Mor}(\Pi(\Sigma),\G(E))$ we
have that $\hat{X}:\pi_1(\Sigma)\rightarrow\G^{X(u_0)}_{X(u_0)}$
is a group homomorphism. If we change $(X,\X)$ by a gauge
transformation $\hat\Psi\in\G^\Sigma$ we get that
$\X_{\hat\Psi}|_{\pi_1(\Sigma)}=\Ad_{\hat\Psi(u_0)}(\X|_{\pi_1(\Sigma)})$,
so that we associate an element in $F(\Sigma,\G(E))$ to the class
in $\M(\Sigma,\G(E))$ represented by $(X,\X)$.

Let us show that this correspondence is injective. Let
$(X_i,\hat{X}_i)$ be two solutions corresponding to the same flat
connection, {\it e.g.}
$\Ad_{\gamma_{21}}(\hat{X}_1|_{\pi_1(\Sigma)})=\hat{X}_2|_{\pi_1(\Sigma)}$,
for some $\gamma_{21}\in\G(E)^{X_2(u_0)}_{X_1(u_0)}$. If we
introduce the local lifting
$\psi_{i\alpha}:U_\alpha\rightarrow\G_{X_i(u_0)}$,
$\psi_{i\alpha}(u)=\hat{X}_i[c^\alpha_{uu_0}]$, it is easy to
verify that
$\hat{\Phi}_\alpha(u)=\psi_{2\alpha}(u)\gamma_{21}\psi_{1\alpha}^{-1}(u)$
for $u\in U_\alpha$ extends to a globally defined map
$\hat{\Phi}\in\G^\Sigma$ such that $s(\hat{\Phi})=X_1$,
$t(\hat{\Phi})=X_2$ and $\X_2=(\X_1)_{\hat\Phi}$.

Let us go in the opposite direction and show that the
correspondence is surjective. In order to do this we first recall
that a $G$-bundle $E\rightarrow B$ is $n$-universal if
$\pi_i(E)=0$ for $i<n$. The following facts are relevant for us:
every $G$-bundle over an $n$-dimensional manifold $N$ is the
pull-back of $E$ for some $X:N\rightarrow B$; moreover if $E$ is
$n+1$-universal, then the $G$-bundles over $N$ are classified by
homotopies from $N$ to $B$ \cite{Steenrod}.

Let now $\rho:\pi_1(\Sigma)\rightarrow\G_{x_0}^{x_0}$ be a flat
connection for some $x_0\in L$; since $\G(E)_{x_0}$ is simply
connected, then $\G_{x_0}\rightarrow L_{x_0}$ is a universal
$2$-bundle for $\G_{x_0}^{x_0}$. This means that the principal
$\G_{x_0}^{x_0}$-bundle $\tilde{\Sigma}\times_\rho\G_{x_0}^{x_0}$
is equivalent to the pull-back $X^*_\rho(\G_{x_0})$ for some
$X_\rho:\Sigma\rightarrow L_{x_0}$. Let
$\Psi_\rho:\tilde{\Sigma}\times_\rho\G_{x_0}^{x_0}\rightarrow\G_{x_0}$
be the bundle map. Finally, define
$\hat{X}_\rho:\pi(\Sigma)\rightarrow\G(E)$ as
$$
\hat{X}_\rho[c_{uv}] =
\Psi_\rho([c^\alpha_{uu_0}],e)\rho[c^\alpha_{u_0u}c_{uv}c^\beta_{vu_0}]\Psi_\rho([c^\beta_{vu_0}],e)^{-1}
~~~~ u\in U_\alpha, v\in U_\beta \;.
$$
It is easy to see that $(X_\rho,\hat{X}_\rho)$ is a well defined
groupoid morphism. \endproof

Remark the crucial role played in the proof by the fact that
$\G_{x_0}$ is simply connected which follows from the assumption that $\G(E)$ is source simply
connected. Indeed the same assumption was used for the Poisson-Lie
group case, see \cite{bonechi}. If the groupoid $\G(E)$ is not
source simply connected or $\dim\Sigma>2$ then the Proposition
\ref{prop_moduli_space_closed} is not true anymore and we only
have the embedding $\M(\Sigma,\G(E)) \subset F(\Sigma,\G(E)) \;.$

\medskip

\begin{rmk}{\rm We can rephrase the above construction by saying that to
every groupoid morphism $(X,\hat{X})$ we can associate a flat
$\G_{x_0}^{x_0}$-bundle over $\Sigma$. In order to be more
explicit, we are going to show that the pull-back principal bundle
$X^*(\G_{x_0})$ admits a flat structure, {\it i.e.} according to
\cite{DeBa} it is isomorphic to
$\tilde{\Sigma}\times_{\hat{X}}\G^{x_0}_{x_0}$. In fact it is
straightforward to verify that $\Phi_X:X^*(\G(E)_{x_0})\rightarrow
\tilde{\Sigma}\times_{\hat{X}}\G_{x_0}^{x_0}$ defined by
$$
\Phi_X(u,\gamma)=([c^\alpha_{uu_0}],\hat{X}[c^\alpha_{u_0u}]\gamma)\;,~~~~
\Phi_X^{-1}([c_{uu_0}],\gamma)=(u,\hat{X}[c_{uu_0}]\gamma) \;
$$
is a well defined principal bundle isomorphism.}\end{rmk}

\medskip

The following examples will help to clarify the above
constructions.

\begin{Exa}\label{MS1}
 Consider Example \ref{FLAT11} and
 let $\G=G$ be a simply connected Lie group seen as a groupoid
over a point $*$; then $\G_*=\G_*^*=G$. Then $G^\Sigma$ can be
seen as a groupoid over a point, hence a group, and
$\bis(G^\Sigma) = G^\Sigma$. The moduli space of solutions
coincides with the moduli space of flat $G$-connections on
$\Sigma$ divided by gauge transformations.
\end{Exa}
\begin{Exa}\label{MS2} There are two extreme cases where the moduli spaces of
solutions is easy to describe. The first one is when $\M$ be a
simply connected symplectic manifold. In fact the (ssc) groupoid
integrating it is $\G(\M)=\M\times \M$ and
$\G_{x_0}^{x_0}=(x_0,x_0)$; the space of flat connection
$F(\Sigma,\G_{x_0}^{x_0})=\{*\}$ is then trivial. The other one is
when $\Sigma=\sphere^2$ where we have that $\M(\sphere^2,\G(E))$
is the space of leaves of $\G(E)$.
\end{Exa}
\begin{Exa}\label{MS3} Let $\G=\M\times G$ be the action groupoid. Then a leaf is an
orbit $L_{x_0}=G/G_{x_0}$, for $x_0\in M$. We have that
$\G^\Sigma=\M^\Sigma\times G^\Sigma$ and thus it is enough to consider $G^\Sigma$
as a gauge group  since the orbits of $\G^\Sigma$ and $G^\Sigma$ coincide. In
fact, the actions of $\hat\Psi=(\psi,\gamma)\in\G^\Sigma$ and
$\gamma\in G^\Sigma$ on $\X$ coincide. This was the gauge group
considered in \cite{bonechi}.
\end{Exa}

\bigskip
\bigskip

\subsection{The open case}\label{moduli_open}

Let $\Sigma$ be a compact surface with boundary. Let us consider
the case with one boundary component $\partial\Sigma$. Let
$\Ep\rightarrow\D$ be a subalgebroid of $E$ and let $\G(\Ep)$ be
the (source connected) subgroupoid of $\G(E)$ integrating it. We
define the relevant moduli space as the space of groupoid
morphisms ${\rm Mor}(\Pi(\Sigma),\G(E); \G(\Ep))$ respecting the
boundary conditions defined in Section \ref{gauge_open} divided by
the action (\ref{gauge_trans_gpd}) of $\G^{\Sigma,\Ep}$, e.g.
\begin{equation}
\label{moduli_space_open} \M(\Sigma,\G(E);\G(\Ep)) = {\rm
Mor}(\Pi(\Sigma),\G(E);\G(\Ep))/\G^{\Sigma,\Ep}\;.
\end{equation}

It is clear that each solution sends the boundary in a fixed leaf
$L\subset\D$ of $\Ep$ and we denote with
$\M(\Sigma,\G(E);\G(\Ep),L)$ the subset of
(\ref{moduli_space_open}) corresponding to this leaf.

In analogy to the closed case, we introduce the following moduli
space of generalized flat connections with the holonomy around the
boundary which takes value in a subgroup. The relevant space is
the union over the leaves $L_x\subset\D$ of $\G(\Ep)$ of the
moduli spaces of flat $\G(E)_x^x$-connections with holonomy around
the boundary living in $\G(\Ep)^x_x$. More precisely, let us
choose $u_0\in\partial\Sigma$ and identify
$\pi_1(\Sigma)=\Pi(\Sigma)_{u_0}^{u_0}$ and denote with
$[\partial\Sigma]\in\pi_1(\Sigma)$ the boundary generator. We
define
\begin{equation}
\label{generalized_flat_open} F(\Sigma,\G(E);\G(\Ep)) =
\bigcup\limits_{x\in\D} \{\rho\in{\rm Mor}(\pi_1(\Sigma),\G(E)^x_x)),\,
\rho[\partial\Sigma]\in\G(\Ep^x_x) \}/ \Ad \;,
\end{equation}
where $\Ad$ is the adjoint action of $\G(\Ep)$.

Let $\{U_\alpha, c^\alpha_{uu_0}\}$ define a trivialization of
$\tilde\Sigma=\Pi(\Sigma)_{u_0}$ as described in Section
\ref{moduli_closed} such that
$c^\alpha_{uu_0}\subset\partial\Sigma$ for all
$u\in\partial\Sigma$.

\medskip
\begin{pro}
\label{proposition_moduli_space_open}
For every surface $\Sigma$ ($\dim \Sigma \leq 2$) with one boundary component and every Lie source simply connected
groupoid $\G(E)$ and Lie source connected subgroupoid $\G(E')\subset\G(E)$ we have that
 $F(\Sigma,\G(E);\G(\Ep)) =
\M(\Sigma,\G(E);\G(\Ep))\;. $
\end{pro}

{\it Proof}. The proof consists in repeating the same steps of the
proof of Proposition \ref{prop_moduli_space_closed} and checking
that the boundary conditions are respected. The map from
$\M(\Sigma,\G(E);\G(\Ep))$ to $F(\Sigma,\G(E);\G(\Ep))$ is defined
in the same way and easily shown to be injective, thanks to the
choice of trivialization of $\tilde\Sigma$. More care is needed
for the inverse map. Let $\rho:\pi_1(\Sigma)\rightarrow\G(E)^x_x$
with $\rho[\partial\Sigma]\in\G(\Ep)^x_x$, for $x\in\D$. We can
then define the $\G(E)^x_x$-bundle
$\tilde\Sigma\times_\rho\G(E)^x_x$ over $\Sigma$ and the
$\G(\Ep)^x_x$-bundle $\partial
\tilde{\Sigma}\times_\rho\G(\Ep)^x_x$ over $\partial\Sigma$. Since
$\G(\Ep)$ is source connected, the $\G(\Ep)^x_x$ bundle
$t:\G(\Ep)_x\rightarrow L_x^{'}$, where $L_x^{'}\subset\D$ is the
leaf containing $x$, is $1$-universal. Then there exists a
$\G(\Ep)^x_x$ bundle map
$\psi_\rho^{'}:\partial\tilde{\Sigma}\times_\rho\G(\Ep)^x_x\rightarrow\G(\Ep)_x$
that can be extended to a $\G(E)^x_x$ bundle map
$\psi_\rho:\tilde\Sigma\times_\rho\G(E)^x_x|_{\partial\Sigma}\rightarrow\G(E)_x$.
Since $\G(E)_x\rightarrow L_x$ is $2$-universal for $\G(E)_x^x$,
we have that $\psi_\rho$ can be extended to the whole bundle over
$\Sigma$ by
$\Psi_\rho:\Sigma\times_\rho\G(E)_x^x\rightarrow\G(E)_x$. The
groupoid morphisms is defined then as in the proof of Proposition
\ref{prop_moduli_space_closed} and respects the boundary
conditions. \endproof

\begin{Exa}
Let us consider $\Sigma=D^1$. Then $\M(\Sigma,\G(E);\G(E'))$ is
the space of leaves of $\G(E')$.
\end{Exa}

We close this Section with a few remarks regarding the moduli
space over the interval $I=[0,1]$. This is closely related to the
explicit construction of the groupoid $\G(E)$ integrating the Lie
algebroid done in \cite{SG} and  \cite{FC}. We consider first the
Lie groupoid morphisms with boundary conditions given by the
trivial groupoid over $\M$ (i.e., $E_0=E_1=\M \times \{ 0\}$ and
$\G(E_0)=\G(E_1)=\M$). Remark that any groupoid morphism $\hat{X}:
I \times I \rightarrow \G(E)$ satisfies these boundary conditions
since they simply mean that $\hat{X}(0,0), \hat{X}(1,1) \in \M$.
The gauge transformations $\hat{\Phi} \in \G^{I,\M, \M}$ are given
by maps $\hat{\Phi}: I \rightarrow \G(E)$ such that
$\hat{\Phi}(0), \hat{\Phi}(1) \in \M$. Then using a standard
argument, it is easy to see that the map
$\hat{X}\rightarrow\hat{X}[1,0]$ defines a bijection between the
moduli space and the groupoid itself, i.e.
$$ {\rm
Mor}( I \times I ,\G(E); \M, \M)/\G^{I,\M, \M} = \G(E).$$

 This description of the
groupoid must be compared with that of \cite{SG,FC}, where the
groupoid is obtained as the space of algebroid morphisms divided
by algebroid homotopies. By taking into account Lemma
\ref{homotopies_gauge_boundary}, it is reasonable to think 
 that $\G^{I,\M, \M}$ is connected to the
identities.

The case of generic boundary conditions can be analogously
treated. Let $\E=\{E_0,E_1\}$ be two subalgebroids of $E$ and let
$\G(E_i)$, $i=0,1$, be the two  subgroupoids of $\G(E)$ integrating
them. It is easy to see that
$$ {\rm
Mor}( I \times I ,\G(E); \G(E_0), \G(E_1))/\G^{I,\E} = \G(E_1)
\backslash t^{-1}(\G(E_1))\cap s^{-1}(\G(E_0))/ \G(E_0).
$$

\bigskip
\bigskip

\section{Concluding remarks}
\label{END}

In this work we have studied the space of Lie algebroid (groupoid)
morphisms modulo gauge transformations. Since our motivations come
from two dimensional topological field theory, the point of view
has been gauge theoretic. In this perspective, we argued that the
choice of finite gauge transformations as the transformations
(\ref{gauge_trans_gpd}) of the groupoid $\G^\Sigma$ is the most
natural one. Indeed the whole story is just a relatively direct
generalization of the group case and the moduli spaces can be
thought of as a generalization of the moduli spaces of flat
connections.

Since groupoids are categories, it is extremely useful to
reconsider the paper from a categorical point of view, where these
choices appear as extremely natural. In fact, let $\H$ and $\G$ be
two groupoids, a groupoid morphism from $\H$ to $\G$ is a
covariant functor. We can consider the functor category ${\mathbf
C}(\H,\G)$, whose objects are the groupoid morphisms ${\rm
Mor}(\H,\G)$ and whose transformations are the natural
transformations between functors; ${\mathbf C}(\H,\G)$ is again a
groupoid. It is easy to verify that the gauge transformations
defined in (\ref{gauge_trans_gpd}) coincide with the natural
transformations. So the moduli space defined in
(\ref{moduli_space_closed}) corresponds to the set $\pi_0({\mathbf
C})$ of the connected components of ${\mathbf C}$ when
$\H=\Pi(\Sigma)$ and $\G=\G(E)$. The content of Proposition
(\ref{moduli_space_closed}) can be expressed by saying that when
$\dim\Sigma\leq 2$ then $\pi_0({\mathbf
C}(\pi_1(\Sigma),\G(E)))=\pi_0({\mathbf C}(\Pi(\Sigma),\G(E)))$.

Moreover, it is important to point out that the in the closed case
the moduli spaces are {\it Morita invariants} of the groupoid
$\G(E)$ (see \cite{BW} for definitions). This fact follows from
the observation that if $\G_1$ and $\G_2$ are Morita equivalent
then also ${\mathbf C}(\pi_1(\Sigma),\G_i)$ are Morita equivalent
groupoids: in particular they have the same space of connected
components.

On a more geometrical side, it will be extremely interesting to
see which geometrical structures can be defined over these moduli
spaces. Indeed some basic facts can be observed now. In Remark
\ref{topology} we pointed out that the moduli spaces are
topological spaces. Moreover, they are the union of moduli spaces
of flat connections and thus they are a collection of symplectic
manifolds (with singularity).

We hope to come back to all these problem in the future and
consider what the quantization of the TFT can bring to the
understanding of these spaces. The TFT which one can associate to
any Lie algebroid is a BF-like theory. In the group case the
quantization of BF-theory gives rise to many interesting
calculations, e.g. the Ray-Singer torsion. It will be interesting
to see if those calculations can be extended to the general case
of Lie groupoids.

%\end{thebibliography}
\end{document}